\documentclass{amsart}%
\usepackage{amsfonts}
\usepackage{amsmath}
\usepackage{amssymb}
\usepackage{graphicx}%
\theoremstyle{plain}

\newtheorem*{proposition}{Main Theorem}

\theoremstyle{definition}
\newtheorem*{acknowledgement}{Acknowledgement}

\theoremstyle{remark}

\numberwithin{equation}{section}
\begin{document}
\title{Estimating the trace-free Ricci tensor in Ricci flow}
\author{Dan Knopf}
\address[Dan Knopf]{ University of Texas at Austin}
\email{danknopf@math.utexas.edu}
\urladdr{http://www.ma.utexas.edu/\symbol{126}danknopf/}
\thanks{Partially supported by NSF grants DMS-0545984 and DMS-0505920.}

\begin{abstract}
An important and natural question in the analysis of  Ricci flow singularity
formation in all dimensions $n\geq4$ is the following: What are the weakest
conditions that provide control of the norm of the Riemann curvature tensor?
In this short note, we show that the trace-free Ricci tensor is controlled in
a precise fashion by the other components of the irreducible decomposition of
the curvature tensor, for all compact solutions in all dimensions $n\geq3$,
without any hypotheses on the initial data.

\end{abstract}
\maketitle

\section{Introduction}

Standard short-time existence results imply that a solution $(\mathcal{M}%
^{n},g(\cdot))$ of Ricci flow on a compact manifold becomes singular at
$T<\infty$ if and only if
\[
\lim_{t\nearrow T}\max_{x\in M^{n}}|\operatorname*{Rm}(x,t)|=\infty.
\]
This suggests the following natural question:\ What are the weakest conditions
that provide control of the norm of the full curvature tensor on a manifold
evolving by Ricci flow?

In any dimension, it is true that a finite-time singularity happens if and
only if
\[
\limsup_{t\nearrow T}\max_{x\in M^{n}}|\operatorname*{Rc}(x,t)|=\infty.
\]
Nata\v{s}a \v{S}e\v{s}um has given a direct proof \cite{Sesum05}. The claim
also follows from independent results of Miles Simon \cite{Simon02} by a short
argument.\footnote{Let $(\mathcal{M}^{n},g(t))$ be a solution of Ricci flow on
a compact manifold such that $g(t)$ is smooth for $t\in\lbrack0,T)$, where
$T<\infty$. If $\limsup_{t\nearrow T}(\max_{x\in\mathcal{M}^{n}}%
|\operatorname*{Rc}(x,t)|)<\infty$, then \cite[Lemma~14.2]{Ham82} guarantees
existence of a complete $C^{0}$ limit metric $g(T)$. One may then apply
\cite[Theorem~1.1]{Simon02}, choosing a background metric $\bar{g}%
:=g(T-\delta)$ such that $(1-\varepsilon)\bar{g}\leq g\leq(1+\varepsilon
)\bar{g}$, where $\varepsilon=\varepsilon(n)$ and $\delta=\delta(\varepsilon
)$. Let $\bar{K}=\max_{x\in\mathcal{M}^{n}}|\operatorname*{Rm}(\bar{g}%
)|_{\bar{g}}$. Simon's theorem implies that there exists $\eta=\eta(n,\bar
{K})$ such that for any $\theta\in\lbrack0,\delta]$, a solution $\hat{g}(s)$
of harmonic-map-coupled Ricci flow exists for $0\leq s<\eta$ and satisfies
$\hat{g}(0)=g(T-\theta)$; moreover, $\hat{g}(s)$ is smooth for $0<s<\eta$.
Since harmonic-map-coupled Ricci flow is equivalent to Ricci flow modulo
diffeomorphisms, the claim follows by taking $\theta=\eta/2$.}

In dimension three, an eminently satisfactory answer is given by the well
known pinching theorem obtained independently by Ivey \cite{Ivey93} and
Hamilton \cite{Ham95}. Their estimate implies in particular that the scalar
curvature dominates the full curvature tensor of any Ricci flow solution on a
compact $3$-manifold with normalized initial data.\footnote{The
Hamilton--Ivery pinching estimate implies the much stronger result that any
rescaled limit of a finite time singularity in dimension three must have
nonnegative sectional curvature.}

Xiuxiong Chen has expressed hope that an appropriate bound on the scalar
curvature might be sufficient to rule out singularity formation in all
dimensions. Partial progress toward this conjecture was made recently by Bing
Wang \cite{Wang07}. He proved that if the Ricci tensor is uniformly bounded
from below on $[0,T)$ and if an integral bound%
\[
\int_{0}^{T}\int_{\mathcal{M}^{n}}\left\vert R\right\vert ^{\alpha}%
\,d\mu\,dt<\infty
\]
holds for some $\alpha\geq(n+2)/2$, then no singularity occurs at time
$T<\infty$.

Recall that in any dimension $n\geq3$, the Riemann curvature tensor admits an
orthogonal decomposition%
\[
\operatorname*{Rm}=U+V+W
\]
into irreducible components%
\[
U=\frac{1}{2n\left(  n-1\right)  }R\left(  g\barwedge g\right)  ,\qquad
V=\frac{1}{n-2}(F\barwedge g),\qquad W=\text{Weyl tensor,}%
\]
where $\barwedge$ denotes the Kulkarni--Nomizu product of symmetric tensors
and $F$ denotes the trace-free Ricci tensor. The purpose of this short note is
to observe that $V$ is always dominated by the other components in the
following sense:

\begin{proposition}
If $(\mathcal{M}^{n},g(\cdot))$ is a solution of Ricci flow on a compact
manifold of dimension $n\geq3$, then there exist constants $c(g_{0})\geq0$,
$C_{1}(n,g_{0})>0$, and $C_{2}(n)>0$ such that for all $t\geq0$ that a
solution exists, one has $R+c>0$ and%
\[
\frac{\left\vert V\right\vert }{R+c}\leq C_{1}+C_{2}\max_{s\in\lbrack
0,t]}\sqrt{\frac{\left\vert W\right\vert _{\max}(s)}{R_{\min}(s)+c}}.
\]

\end{proposition}

\section{Proof of the main theorem}

Define%
\[
a=\left\vert F\right\vert =\frac{\sqrt{n-2}}{2}\left\vert V\right\vert ,
\]
noting that $a$ is smooth wherever it is strictly positive. Choose $c\geq0$
large enough so that $R_{\min}(0)+c>0$ and define%
\[
b=R+c,
\]
noting that $b>0$ for as long as a solution exists.

In any dimension $n\geq3$, one has%
\[
\frac{\partial}{\partial t}\left\vert F\right\vert ^{2}=\Delta\left\vert
F\right\vert ^{2}-2\left\vert \nabla F\right\vert ^{2}+\frac{4(n-2)}%
{n(n-1)}R\left\vert F\right\vert ^{2}-\frac{8}{n-2}\operatorname*{tr}%
F^{3}+4W(F,F),
\]
where $\operatorname*{tr}F^{3}=F_{i}^{j}F_{j}^{k}F_{k}^{i}$ and
$W(F,F)=W_{ijk\ell}F^{i\ell}F^{jk}$. It follows from Cauchy--Schwarz that $a$
obeys the differential inequality%
\[
a_{t}\leq\Delta a+\frac{2(n-2)}{n(n-1)}a(b-c)-\frac{4}{n-2}a^{-1}%
\operatorname*{tr}F^{3}+2a^{-1}W(F,F).
\]
The positive quantity $b$ evolves by%
\[
b_{t}=\Delta b+2a^{2}+\frac{2}{n}(b-c)^{2}.
\]

To prove the theorem, it will suffice to bound the scale-invariant
non-negative quantity%
\[
\varphi=\frac{a}{b}.
\]
Because $\Delta\varphi=b^{-1}(\Delta a-\varphi\Delta b)-2\left\langle
\nabla\varphi,\nabla\log b\right\rangle $, one has%
\[
\varphi_{t}\leq\Delta\varphi+2\left\langle \nabla\varphi,\nabla\log
b\right\rangle +2\rho\varphi,
\]
where the reaction term is%
\[
\rho=\frac{n-2}{n(n-1)}(b-c)-\frac{2}{n-2}\frac{\operatorname*{tr}F^{3}}%
{a^{2}}+\frac{W(F,F)}{a^{2}}-\frac{(b-c)^{2}}{nb}-a\varphi.
\]
There exist positive constants $c_{1},c_{2}$ depending only on $n\geq3$ such
that%
\[
\left\vert \frac{2}{n-2}\operatorname*{tr}F^{3}\right\vert \leq c_{1}%
a^{3}\qquad\text{and}\qquad\left\vert W(F,F)\right\vert \leq c_{2}\left\vert
W\right\vert a^{2}.
\]
Hence%
\[
\rho\leq\frac{n-2}{n(n-1)}(b-c)+c_{1}a+c_{2}\left\vert W\right\vert
-\frac{(b-c)^{2}}{nb}-a\varphi.
\]

Define constants $\alpha,\beta,\gamma$ by $\alpha^{2}=c_{1}$, $\beta^{2}%
=\frac{n-2}{n(n-1)}$, and $\gamma^{2}=c_{2}$. Fix $\varepsilon>0$ and choose
$C_{1}=\max\{\alpha^{2}+\beta,\varphi_{\max}(0)+\varepsilon\}$ and
$C_{2}=\gamma$. Consider the barrier function%
\[
\Phi(t)=C_{1}+C_{2}\max_{s\in\lbrack0,t]}\sqrt{\frac{\left\vert W\right\vert
_{\max}(s)}{b_{\min}(s)}},
\]
noting that $\Phi$ is monotone nondecreasing. If $\varphi_{\max}(t)\geq
\Phi(t)$ at some $t>0$, then at any $(x,t)$ where $\varphi$ attains its
spatial maximum, one has%
\[
a\geq(\alpha^{2}+\beta)b+\gamma\sqrt{b\left\vert W\right\vert },
\]
which implies that
\begin{align*}
a^{2}  &  \geq\alpha^{2}ab+\beta^{2}b^{2}+\gamma^{2}b\left\vert W\right\vert
\\
&  \geq\frac{n-2}{n(n-1)}(b^{2}-bc)+\left(  c_{1}a+c_{2}\left\vert
W\right\vert \right)  b-\frac{1}{n}(b-c)^{2},
\end{align*}
hence that $\rho\leq0$, hence that $\varphi_{t}\leq0$, hence that $\frac
{d^{+}}{dt}\varphi_{\max}(t)\leq0$, understood in the usual sense as the
$\limsup$ of difference quotients. It follows that $\varphi_{\max}(t)\leq
\Phi(t)$ for as long as a solution exists.

\begin{acknowledgement}
The author warmly thanks Matt Gursky for a number of stimulating discussions.
\end{acknowledgement}

\end{document}